\documentclass{aupl}

\usepackage{
amsfonts,
amsmath,
latexsym,
amssymb,
enumerate,
verbatim,
mathrsfs,
graphicx,
centernot,
accents, 
color,
}

\usepackage[all]{xy}

\newcommand{\labbel}[1]{\label{#1} [[{\bf #1}]]}  
\renewcommand{\labbel}{\label}

\newcommand{\red}[1]{{#1}}

\newcommand{\contact}{{contact$^{ \text{-}}$} } 
\newcommand{\Contact}{{Contact$^{ \text{-}}$} }

\newcommand{\ddd}{\mathrel{\delta}} 
\newcommand{\nd}{\mathrel{\centernot\delta}}

\newcommand{\nup}{{\centernot\uparrow}}
\newcommand{\da}{{\downarrow}}
\newcommand{\ua}{{\uparrow}}

\newtheorem{theorem}{Theorem}

\newtheorem{proposition}[theorem]{Proposition} 
 
\newtheorem{corollary}[theorem]{Corollary}

\newtheorem*{claim*}{Claim}

\newtheorem*{theorem*}{Theorem}
\newtheorem*{proposition*}{Proposition}
\newtheorem*{corollary*}{Corollary}
\newtheorem*{lemma*}{Lemma}
\newtheorem*{scholion*}{Scholion}

\theoremstyle{definition}
\newtheorem{definition}[theorem]{Definition}

\theoremstyle{remark}
\newtheorem{remark}[theorem]{Remark}

\newtheorem*{remark*}{Remark}
\newtheorem*{remarks*}{Remarks}

\newtheorem*{observation*}{Observation}

 \allowdisplaybreaks[1]

\numberwithin{equation}{section}

\begin{document}

 \title{Contact posets}

\author{Paolo Lipparini} 
\address{Dipartimento di Matematica\\Viale 
della  Prossima Pensione
\\Universit\`a di Roma ``Tor Vergata'' 
\\I-00133 ROME ITALY}

\email{lipparin@axp.mat.uniroma2.it}


\begin{abstract}
We study contact posets and show
that every contact poset can be embedded into a 
Boolean poset with overlap contact relation.
Contact posets and (nonadditive) contact
semilattices have the superamalgamation property,
Fra\"\i ss\'e limits and model completion.
\red{Some results apply to event structures
with binary conflict, as introduced in computer science.} % fine red 
\end{abstract} 


\maketitle

Proximity structures provide a useful generalization of topology
\cite{D,N}. Their algebraic analogue has been 
considered under various names, more frequently, \emph{contact algebras}.
Contact (Boolean) algebras play an important role in 
 region-based theory of space \cite{APB}.
Algebras with less structure, such as 
contact lattices, have been studied; more generally,
\cite{I} provides arguments suggesting the naturalness 
of contact join-semilattices.
Another argument suggesting that it is natural to consider
only the partial order or the join operation is that 
meet and complements are not preserved 
by the image function from $\mathcal P(X)$ 
to $\mathcal P(Y)$ associated to some function
from $X$ to $Y$; see \cite{mtt}
for further elaborations.  
More references and details can be found in the quoted sources.

Here we go one step further and study \emph{contact posets}, proving some
representation theorems. 
\red{Under the name \emph{weak contact structure},
the definition of a contact bounded poset appeared in 
\cite{DW}, but there the authors studied only the
case of distributive lattices.
A similar notion, with an additional finitary condition,
appeared earlier in the theory of concurrent systems in computer science \cite{WN}.
See Remark \ref{eve}  below for more details.} % fine red 

The axiom of choice is not used in the present note.

\begin{definition} \labbel{def}    
By a \emph{poset} we mean a partially ordered set
with a minimum element $0$.
The existence of $0$ is assumed only in order to simplify
the statements of the representation theorems, 
see Remark \ref{senzazero} below. 
A \emph{\contact relation on a poset} $\mathbf P= (P, {\leq})$
is a binary relation $\delta$  on $P$ such that     

 \begin{align}
\labbel{sym}    \tag{Sym} 
 &  a \ddd b \Leftrightarrow    b \ddd a
\\
\labbel{emp} \tag{Emp}
& a \ddd b  \Rightarrow a > 0  \  \& \ b>0,  
\\
\labbel{ext}    \tag{Ext} 
& a \ddd b \  \& \  a \leq a_1\  \&\  b \leq b_1
\Rightarrow a_1 \ddd b_1,
\\
\labbel{ref}    \tag{Ref} 
& n \neq 0  \Rightarrow  n \ddd n,
  \end{align}
for all $n, a, b, a_1, b_1 \in P$.
In the literature sometimes the name \emph{weak contact}
is used in place of contact$^{ \text{-}}$. 
We write $a \nd b$ to mean that  $a \ddd b$ does not
hold.  

A \emph{\contact  poset}
is a structure   $(P, {\leq}, { \delta  } )$,
where $(P, {\leq})$ is a poset and $\delta$ 
is a \contact relation on $(P, {\leq})$;
\emph{\contact semilattices}, \emph{\contact Boolean algebras},
etc., are defined in an analogous way. 
The reason for the minus sign
in the word \contact  shall be explained
soon. 
Semilattices shall always be considered as \emph{join semilattices},
namely, the associated order is given by $a \leq b$ if $a+b=b$.  

From \eqref{ext} and \eqref{ref}  we get
  \begin{align}
\labbel{inh}    \tag{Inh} 
  n \neq 0  \  \&\   n \leq a\  \&\  n \leq b
\Rightarrow 
 a \ddd b. 
  \end{align}

Notice that if  $a \nd b$, then the meet
$ab$ of $a$ and $b$ exists in $\mathbf P$, in fact, $ab=0$,
because of \eqref{ref} and \eqref{ext}.     
If $\mathbf P$ is a poset with $0$,
then, setting $a \ddd b$ if there is $n \in P$, $n >0$
such that $n \leq a$ and $n \leq b$,
we get a \contact relation, which is called the
\emph{overlap} (or \emph{minimal}
or \emph{trivial}) \contact relation on  $\mathbf P$. 

Most authors include the following 
additivity condition in the definition
of a contact structure:
  \begin{align}
\labbel{add}    \tag{Add} 
  a \ddd b+c 
\Rightarrow 
 a \ddd b \text{ or } a \ddd c.
  \end{align}

We shall generally not include \eqref{add} in the definition 
of \contact structures, since \eqref{add} is not
even expressible in the setting of posets.  
Moreover, there are structures of possible interest in which
\eqref{add} is expressible, but not satisfied; 
just to mention an elementary example, \eqref{add}
fails in  the $5$-element
modular lattice $\mathbf M_3$ with $3$ atoms $a$, $b$, $c$, with the overlap  
relation.
The additivity condition \eqref{add}
fails  in $\mathbf M_3$ even when 
$a \ddd b$ is added as the only nontrivial pair in contact.
 As well-known, there are many interesting
examples of nondistributive lattices, sometimes with applications 
outside mathematics \cite{PP,R}. 

When we assume that some contact
semilattice is additive, we shall explicitly
include the word additive.
In the general case, we use a minus sign in the word \contact in order to remind
the reader of our convention.

An \emph{embedding} $\varphi$  is an injective order-preserving 
and $0$-preserving map
such that $a \ddd b$ if and only if $ \varphi (a) \ddd \varphi(b) $,
for all elements $a$ and $b$ in the domain. 
When dealing with semilattices, $\varphi$  is also assumed to preserve 
$+$.    

We now prove some representation theorems.
\end{definition}

\begin{proposition} \labbel{emb1}
Every \contact poset (semilattice)
can be embedded into a poset (semilattice)
with overlap \contact relation.
 \end{proposition}

  \begin{proof}
If $(P, {\leq}, { \delta  } )$ is a \contact poset,
consider the function $\varphi: P \to \mathcal P(P)$  
defined by $\varphi(a) = \nup a = \{ \, x \in P  \mid
  a \centernot \leq x\, \} $.
Consider the \contact poset  $\mathbf Q = (Q, { \subseteq }, { \delta  }_Q )$,
where $Q= Im \varphi \cup \{ \, \nup a \cap \nup b   \mid 
a,b \in P, a \ddd b \, \} $ and, for $x,y \in Q$, 
$x \ddd _Q y$ if  there is $q \neq \emptyset $,
$q \in Q$ such that    $ q \subseteq x $ and $ q \subseteq   y$.
It is elementary to see that $\mathbf Q$ is a
\contact poset; and  $\delta_Q$ is  the overlap relation
by construction.
Notice that $ \emptyset = \varphi (0)$ is the minimum of    $\mathbf Q$.

Now consider $\varphi$  as a function from $P$ to $Q$.
We claim that $\varphi$  is an embedding. 
It is standard to see that $\varphi$  is an order-embedding.
If $a,b \in P $ and $ a \ddd b$, then $a > 0$
and $b > 0$, by    \eqref{emp}, thus  
$\nup a \cap \nup b \neq \emptyset  $, since
$0 \in \nup a \cap \nup b$. Since
$ a \ddd b$, then $\nup a \cap \nup b  \in Q$,
by the definition of $Q$,  
hence $\varphi(a) \ddd _Q \varphi (b)$,
since $\varphi(a) = \nup a \supseteq  \nup a \cap \nup b $
and similarly  $\varphi(b) \supseteq  \nup a \cap \nup b $,
so we can take $q= \nup a \cap \nup b$.

Conversely,  suppose that 
$\varphi(a) \ddd _Q \varphi (b)$, thus 
$q \subseteq \varphi(a) = \nup a$ and 
$q \subseteq \varphi(b) = \nup b$, 
for some $q \in Q$, $q \neq \emptyset $,
by the definition of $\ddd _Q$.  There are two cases.
(i) If $q \in Im \varphi $, then 
$q = \nup c$, for some $c \in P$.
From    $\nup c = q \subseteq \nup a$,
by taking complements, 
we get $\ua c \supseteq \ua a$,
where $\ua a = \{ \, x \in P  \mid
 x \geq a \, \} $. 
 Since
$a \in \ua a$, then $a \in \ua c$, that is,     
$c \leq a$. Similarly, $c \leq b$,
thus $a \ddd b$, by \eqref{inh}, since $c > 0$
(indeed, if $c=0$, then $q= \nup c = \emptyset  $,
contrary to our assumption).
We now deal with the other case, that is,
(ii)  $q = \nup c \cap \nup d$,
for $c,d \in P$ such that $c \ddd d$.   
From 
$ \nup c \cap \nup d = q \subseteq \varphi(a) = \nup a$,
again taking complements, we get 
$ \ua c \cup \ua d  \supseteq  \ua a$,
hence either $a \in \ua c $ or  $a \in \ua d $, that is,
either $c \leq a$ or $d \leq a$.
Symmetrically, either 
$c \leq b$ or $d \leq b$.
Since $c \ddd d$, then $c >0$ and  $d >0$,
by \eqref{emp}.   
Hence, if both $c \leq a$ and  $c \leq b$,
then $a \ddd b$, by  \eqref{inh}.
Similarly,  $a \ddd b$ if both $d \leq a$ and  $d \leq b$.
 Otherwise, say, 
both $c \leq a$ and  $d \leq b$,
then  $a \ddd b$, by \eqref{ext} and  
$c \ddd d$.
The remaining case 
$d \leq a$ and  $c \leq b$ is similar.

Now suppose that  $(P, {+}, { \delta  } )$ is a \contact semilattice.
Define $\varphi$  and $Q$ as above.
Here, $Q$ is not necessarily closed under union, hence
we need consider $Q^+$, the closure of $Q$ under finite unions.
Namely, $Q^+$ is the set of subsets of $P$ having the form
$q_1 \cup \dots \cup q_n$, for  $q_1 , \dots,  q_n \in Q$.
As above, define 
$x \ddd _{Q^+} y$ if  there is $r \neq \emptyset $,
$r \in Q^+$ such that    $ r \subseteq x $ and $ r \subseteq   y$
and consider the \contact semilattice 
$(Q^+, {\cup}, { \ddd _{Q^+}  } )$.
Again, $\ddd _{Q^+}$ is the overlap  
relation by construction. 
The same proof as above shows that 
$\varphi$, considered as a function from 
$P$ to $Q^+$, is an embedding of 
  \contact semilattices.
Indeed, if $ r=q_1 \cup \dots \cup q_n \in Q^+$
and $ r \neq \emptyset  $, then at least one $q_i$ is not $ \emptyset $.   
If both $r \subseteq \varphi (a)$
and  $r \subseteq \varphi (b)$, then
both 
$q_i \subseteq \varphi (a)$ and
$q_i \subseteq \varphi (b)$,
hence we can repeat the above argument
with $q_i$ in place of $q$.  
 \end{proof}

The proof of Proposition \ref{emb1}
shows a little more. 

\begin{corollary} \labbel{cor}
Every \contact poset $\mathbf P$  
can be embedded into 
(the order-reduct of a)
semilattice
with overlap \contact relation,
in such a way that all the existing, possibly infinitary,
joins in $\mathbf P$ are preserved. 
 \end{corollary} 

\begin{theorem} \labbel{thm}
(a) Any \contact poset can be embedded into 
a \contact poset which is the order-reduct of a 
complete atomic Boolean algebra with overlap contact relation.

(b) Any \contact semilattice can be (semilattice-) embedded into 
a bounded complete \contact lattice with overlap \contact relation.
 \end{theorem} 

\begin{proof}
(a) In view of Proposition \ref{emb1}, 
and since the composition of two embeddings is an embedding, 
it is enough to show that a \contact poset with overlap \contact relation
can be embedded into a  contact Boolean poset satisfying the
further requested properties.
So let $(Q, {\leq}, { \delta  } )$ be a \contact poset
with overlap \contact relation.
Let $R= \mathcal P(Q \setminus \{ 0 \})$ and 
$\mathbf R = (R, { \subseteq }, \ddd _R)$,
where  $\ddd _R$ is the overlap relation,
given by $x \ddd _R y$ if $x \cap y \neq \emptyset $.   
Let $\psi : Q \to R$ be defined by
$\psi (q) = \da q  = \{ \, x \in Q  \mid
0 \neq x \leq q \, \} $.   

The only nontrivial thing to check is that 
$\psi$ is a $\delta$-embedding.  
If $a,b \in Q$ and $a \ddd_Q b$,
then there is $c \in Q$, $c >0$  such that $c \leq a$
and $c \leq b$, since    $\ddd_Q$ is the overlap relation.
Then $\psi (c) \neq \emptyset  $,
$\psi (c) \subseteq \psi (a)$ 
and $\psi (c) \subseteq \psi (b)$,
thus $\psi (a) \ddd_R \psi (b)$, by \eqref{inh}.  
Conversely, if 
 $a \nd_Q b$, then necessarily $0$ is the meet of $a$ and  $b$,
again by \eqref{inh}. Since $\psi$ preserves existing meets,
then $\psi (a) \cap \psi (b) = \emptyset $,
thus  $\psi (a) \nd_Q \psi (b)$, since
$\nd_Q$ is the overlap relation.  
  
(b) As well-known, every poset $\mathbf Q$  can be embedded into a 
bounded complete
lattice by an embedding $\chi$ which preserves all 
meets and joins existing in $\mathbf Q$ \cite{H}.
It follows that if, furthermore, $\mathbf Q$ is a join semilattice,
then $\chi$ is a semilattice embedding.
Again, in view of Proposition \ref{emb1}, 
it is no loss of generality to assume that  
$\mathbf Q$ has overlap \contact relation.
Then the fact that 
$\chi$ is an order-embedding which preserves existing meets
is enough to show that the arguments in the last paragraph 
of the proof of (a) work.
 \end{proof}

We do not know whether results analogue to 
Proposition \ref{emb1} and Theorem \ref{thm}(b)
can be proved for additive contact semilattices.  
For sure, the items
in Theorem  \ref{thm} cannot be joined together 
in order to show that every \contact semilattice can be embedded into 
(the semilattice reduct) of a distributive lattice
with overlap \contact relation.
In fact, 
 if some \contact semilattice $\mathbf S$  can be embedded into 
 a distributive overlap lattice,
then $\mathbf S$ is additive.
Indeed, if $ a \ddd b+c$ in some distributive lattice
 with overlap \contact relation,
then $0 < a(b+c)= ab +ac$, so that either
$0 < ab$ or $0< ac$, hence   
 either $ a \ddd b$ or $ a \ddd c$. 
Properties stronger than additivity are needed in order
to get embeddability in a distributive lattice with
overlap contact.
Characterizations of contact semilattices embeddable into 
 a contact distributive lattice are presented in \cite{I,cs}. 

Recall that a class $\mathcal K$ of structures of the same type and
closed under isomorphism has the \emph{strong amalgamation property}  
if, whenever $\mathbf A,\mathbf  B, \mathbf  C \in  \mathcal K$, 
 $ \mathbf  C \subseteq \mathbf A $, 
$ \mathbf  C \subseteq \mathbf B $ and $A \cap B = C$,
then there is some structure $\mathbf  D \in \mathcal K$ 
such that $\mathbf A \subseteq \mathbf  D$
and $\mathbf  B \subseteq \mathbf  D$.  
If $\mathcal K$ is a class  of ordered structures,
the \emph{superamalgamation property} means that we can
also obtain that, if $a \in A$, $b \in B$ and
$a \leq _{ \mathbf  D} b $, then there is $c \in C$ 
such that   $a \leq _{ \mathbf  A} c \leq _{ \mathbf  B} b $,
and symmetrically when $b \leq _{ \mathbf  D} a $.
For contact (additive Boolean) algebras, the strong amalgamation property (under the
terminology \emph{disjoint amalgamation property})
is proved in \cite{DL}.

We refer to \cite[Section 7.1]{H}
for the notions of a Fra\"\i ss\'e limit and of a model completion.
As a way of example, we just mention that the random graph is the Fra\"\i ss\'e 
limit of the class of finite graphs.

\begin{theorem} \labbel{ap}
The theories of \contact posets and of
\contact semilattices have the superamalgamation property. 
In each case, the class of finite models has
a Fra\"\i ss\'e limit $\mathbf M$.
In each case, the first-order theory of $\mathbf M$
 is 
$ \omega$-categorical, has quantifier elimination and
is the model completion of the respective theory.
 \end{theorem} 

\begin{proof} 
Given a triple $\mathbf A$, $\mathbf  B$, $\mathbf  C$ 
of \contact posets to be  amalgamated, let $\leq$ be
the smallest relation on $D= A \cup B$
such that $\leq$ extends the orders on $\mathbf A$
and $\mathbf  B$ and also contains 
${\leq_ {\mathbf  A}} \circ {\leq_ {\mathbf  B}}$
and ${\leq_ {\mathbf  B}} \circ {\leq_ {\mathbf  A}}$.
In \cite[Lemma 3.3]{Jo} it is shown that 
$\leq$ is a partial order on $D$ and that $(D, \leq_ {\mathbf  D})$
strongly amalgamates the order-reducts of $\mathbf A$ and $\mathbf  B$
over $\mathbf  C$. Superamalgamation
holds by the very definition of $\leq$. 
 
Now define $\delta$ on $D$ by
$d \ddd e$ if there are $a,b \in D$
such that $ a \leq d$, $b \leq e$ and
either (a) $a,b \in A$ and $a \ddd _{ \mathbf A}  b$,
or  (b) $a,b \in B$ and $a \ddd _{ \mathbf B} b$.
It is easily verified that $\mathbf  D$ becomes a 
\contact poset. Moreover, if, say, $d,e  \in A$,
then $d \ddd e$ if and only if $d \ddd _{ \mathbf A}  e$.
Sufficiency is immediate from the definition of $\delta$ in $\mathbf  D$;
on the other hand, if $d,e  \in A$,
and $d \ddd e$ is witnessed by 
$ a \leq d$ and $b \leq e$ with, say, $a,b \in B$
such that $a \ddd _{ \mathbf B} b$,
then  there are $c,c_1 \in C$ such that  
$ a \leq_{ \mathbf B} c \leq _{ \mathbf A} d$ and
 $b _{ \mathbf B}\leq c_1 \leq  _{ \mathbf A}  e$,
by the definition of $\leq$ in $\mathbf  D$. 
From $a \ddd _{ \mathbf B} b$ and 
$ a \leq_{ \mathbf B} c $,
 $b _{ \mathbf B}\leq c_1$,
we get $c \ddd _{ \mathbf B} c_1$,
by applying \eqref{ext} in $\mathbf  B$.
Since $\mathbf  C$ embeds in $\mathbf  B$,
then $c \ddd _{ \mathbf C} c_1$,
hence $c \ddd _{ \mathbf A} c_1$,
since $\mathbf  C$ embeds in $\mathbf A$.
Then, working in $\mathbf A$, we get 
$d \ddd _{ \mathbf A}  e$ from
$ c \leq _{ \mathbf A} d$,
 $c_1  \leq_{ \mathbf A}  e$
and \eqref{ext}. The other cases are similar
or easier. We have showed that the inclusions are embeddings
from $\mathbf A$, respectively, $\mathbf  B$ to $\mathbf  D$,
thus $\mathbf  D$ is a superamalgamating structure
in the case of \contact posets.
 
In the case of \contact semilattices,
proceed as above with regard to the order and
the \contact structure. Of course, $\mathbf  D$ is not necessarily a
semilattice, but \cite[p.\ 205]{Jo} proves that  joins of $\mathbf A$ are preserved in
$(D, \leq_ {\mathbf  D})$, and similarly for  joins of $\mathbf B$.
By Corollary \ref{cor}, $\mathbf  D$ can be order- and contact-embedded
into  some \contact semilattice $\mathbf  E$ in such a way that existing
joins in $\mathbf  D$ are preserved, thus joins of $\mathbf A$
are preserved in going from $\mathbf A$ to $\mathbf  E$  
 and similarly for
$\mathbf  B$. This means that $\mathbf  E$  
superamalgamates $\mathbf A$ and $\mathbf  B$ over $\mathbf  C$ 
in the class of \contact semilattices.

In both cases, the class of finite structures has
a Fra\"\i ss\'e limit by
\cite[Theorem 7.1.2]{H}. 
Since both classes are locally finite, the remaining parts
of the theorem follow
by standard arguments, e.~g., \cite[Theorem 7.4.1]{H}.
\end{proof}

\Contact distributive lattices (no matter whether additive or not)
fail to have the amalgamation property (see
\cite{KMPT} for details about the various 
versions of the amalgamation property). 
Indeed, let $\mathbf  C$ be a $3$-element chain
with elements $0, c, 1$, let both $\mathbf A$ 
and $\mathbf  B$ be obtained by adding a complement 
($a$, $b$, respectively) to $c$.
Let $\delta$ be the overlap relation in $\mathbf A$, thus 
$a \nd c$, while  let  $b \ddd c$ in $\mathbf  B$.    
In any amalgamating distributive lattice
$a$ and  $b$ should be identified,
since complements are unique in distributive lattices, but this contradicts 
the embedding property, since $a$ and $b$
satisfy distinct relations.  

The above argument is really general:
as far as we have a class $\mathcal K$ of distributive lattices
in an expanded language such that  (a) there is some $\mathbf C \in \mathcal K$
and $c \in L$ without complement in $\mathbf  C$,
(b) $\mathbf  C$ has some extension
$\mathbf A \in \mathcal K$
in which $c $ has a complement, and
(c) the extra (non-lattice) structure of $\mathbf A$
is not uniquely determined by the extra structure of $\mathbf  C$,
then $\mathcal K$ has not the amalgamation property.   
What lies behind the argument is the fact that 
distributive lattices have the amalgamation property 
but not the strong amalgamation property.

\begin{remark} \labbel{senzazero}
Axioms can be given for \contact posets and semilattices
without $0$. 

A \emph{\contact  poset (semilattice) without $0$} 
is a poset (semilattice) with
 a binary relation $\delta$  such that     
\eqref{sym},  \eqref{ext} and
 \begin{align}
\labbel{ref*}    \tag{Ref$^*$} 
&  n \ddd n
  \end{align}
 hold, for every element $n$.

The definition is actually simpler than
Definition \ref{def}; on the other hand,
it is simpler to deal with the $0$-ed notions
when we are concerned  with embeddings into,
say, Boolean algebras.  
 \end{remark}   

\red{
\begin{remark} \labbel{eve}
Event structures have been studied in computer
science in various distinct formulations; 
 see \cite{GP}
for a discussion.
In the terminology from 
\cite[Section 8]{WN} an
\emph{event structure}  $(E, \leq, \#)$ 
is a poset together with a binary symmetric irreflexive relation $\#$  such that
  \begin{enumerate}    
\item
(down-finiteness)
 for every $e \in E$ the set of all the $\leq$-predecessors of $e$ is finite, and
\item
 $e \mathrel{\#} e'$ and $e' \leq e''$ imply    $e \mathrel{\#} e''$.
   \end{enumerate} 
In more recent terminology, event structures in the above sense
are frequently named \emph{prime event structures with binary conflict}.

Considering the dual order $\leq'$ and letting $\delta$ be the negation of 
$\#$, that is  $e \ddd e'$ if and only if  not $e \mathrel{\#} e'$,
an event structure is exactly a 
\contact  poset  without $0$, as defined in the above remark,
with the further property that up-finiteness holds, that is,
the set of all the $\leq'$-successors of each element is finite.
Just notice that the contrapositive of (2) reads:
 $e \ddd  e''$ and $e'' \leq' e'$ imply    $e \ddd e'$.

In the proof of Theorem \ref{ap} for \contact posets
we have taken $D=A \cup B$, hence if up-finiteness
holds both in $\mathbf A$ and $\mathbf  B$, then it holds in 
$\mathbf  D$. Thus we get a proof for the following theorem.
  \end{remark}   

\begin{theorem} \labbel{apevent}
The class of event structures in the sense of Definition \ref{eve},
with respect to the standard model theoretical notion of embedding,
has the strong amalgamation property, more generally, the superamalgamation property. 
 \end{theorem}   

On the other hand, the construction of a Fra\"\i ss\'e limit
uses infinite unions, hence up-finiteness is not preserved. 
Let us also remark that authors in computer science
generally use a different notion of  morphism, see e.~g., 
\cite[Subsection 8.1]{WN}. With an appropriate notion 
of embedding, universal homogeneous models generally
exist
 \cite{Dr}.
} % fine red

\begin{remark} \labbel{inpiu}
Our main motivation for the study of \contact posets
comes from some more complex aspects.
Consider a poset $\mathbf P$ with $0$  endowed with
a closure (= isotone, idempotent and extensive) operation $K$, a typical example
is $(\mathcal P(X), {\subseteq} , K)$,
where $X$ is a topological space with closure operation $K$.  
One can define a \contact relation $\delta$ on $P$  
by setting $a \ddd b$ if  there exists $n \in P$
such that $n > 0$,  $n \leq Ka$ and $n \leq Kb$.   
In the case of topological spaces this is called 
the \emph{standard proximity} \cite[Example 2.1.3]{D}.
Notice that $\delta$ induces an overlap relation
on the set of closed elements of $\mathbf P$
($c$ is \emph{closed} if $Kc=c$).  Compare
also the example in the next remark.
  
The study of \contact posets and of overlap relations
is thus preparatory for the study
of structures of the form 
$(P, {\leq}, K, { \delta } )$
and their representations,
an endeavor we are currently working on and
which is motivated by \cite{mtt}.
 \end{remark}

\emph{Concluding remark.}
In view of Proposition \ref{emb1}
and Theorem \ref{thm}, one might argue that the theory of
\contact semilattices is too weak to really distinguish between
overlap and nonoverlap contact relations.
This is certainly a reasonable point of view.

On the other hand, the subject might be seen from
a different perspective.
For example, think of the metric proximity
$\delta$ defined over the metric space of the rationals in $[0,1]$,
with the usual euclidean metric.
The proximity $\delta$ is defined for subsets of $X= \mathbb Q \cap [0,1]$  
  by $a \ddd b$ if 
$ \inf \{ \, d(x,y) \mid x \in a, y \in b \,\} = 0$ 
  \cite[Example 2.1.4]{D}. In algebraic terms,
 $\mathcal P(X)$ becomes endowed with the  structure of a
contact (Boolean) algebra.

There is a natural homomorphism from
the contact semilattice  $\mathbf S
= (\mathcal P(X), \allowbreak {\cup}, {\ddd})$
to a semilattice with overlap relation.
Let $\mathbf T
= (T, {\cup}, {\ddd_T})$,
where $T$ is the set of closed
subsets of $[0,1]$, the unit interval of real  numbers, 
and $ \delta _T $ is the overlap relation.
Compare the previous example.
Then the function $\varphi$  which sends
a subset $a $ of  
$\mathbb Q \cap [0,1]$ to the closure of $a$
in $[0,1]$ is a homomorphism from
$\mathbf S$ to $\mathbf T$.
   
In this sense, Theorem \ref{thm}(a) shows that,
given an arbitrary \contact poset $\mathbf P$, 
one can always add ``imaginary'' elements to $P$
 in such a way that the contact between two proximal
elements can be always witnessed by an actual
``point of contact''---or, better,
region of contact---which is contained in both elements.  

From an epistemological point of view,
the above example supports the thesis that 
real numbers, rather than sharing some kind of ``real''
existence,  are nothing more than ``ideal objects''
added in order to obtain suitable completions.
In this sense, under many respects, metrics or proximities over the rational numbers
are a much more concrete version of (and, frequently, equivalent to)  
the more usual topology on $\mathbb R$.

\end{document}